# A multiresolution triangular plate-bending element method

YiMing Xia

(Civil Engineering Department, Nanjing University of Aeronautics and Astronautics, Nanjing, P.R.China)

email:xym4603@sina.com

**Abstract**

A triangular plate-bending element with a new multi-resolution analysis (MRA) is proposed and a novel multiresolution element method is hence presented. The MRA framework is formulated out of a displacement subspace sequence whose basis functions are built out of scaling and shifting on the element domain of basic full node shape function. The basic full node shape function is constructed by means of extending the shape function triangle domain for a split node at the zero coordinates to the hexagon area enclosing the zero coordinates. As a result, a new split-full node notion is presented and a novel rational MRA concept together with the resolution level (RL) is constituted for the element. Via practical examples, it is found that the traditional triangular plate element and method is a mono-resolution one and also a special case of the proposed element and method. The meshing for the monoresolution plate element model is based on the empiricism while the RL adjusting for the multiresolution is laid on the rigorous mathematical basis. The analysis clarity of a plate structure is actually determined by the RL, not by the mesh. Thus, the accuracy of a structural analysis is replaced by the clarity, the irrational MRA by the rational and the mesh model by the RL that is the discretized model by the integrated. The continuous full node shape function unveils secrets behind assembling artificially of node-related items in global matrix formation by the conventional FEM.

## 1. Introduction

Multi-resolution analysis (MRA) is a popular technique that has been applied in many domains such as the signal and image processing, the damage detection and health monitoring, the differential equation solution, etc. Herein, MRA is referred to a method by which the amount of exposed details (nodes) over a concerned area can be modulated freely in a certain manner at a request. Meanwhile, finite element method (FEM) [1] is also a robust numerical analysis tool that has been widely applied in the differential equation solution, the engineering structural analysis, the thermodynamics etc. In fact, FEM has always employed the MRA technique to deal with practical problems indiscernibly. As what is commonly known, when encountering a practical engineering problem, the selection of an appropriate meshing scheme (a node distribution pattern) will be a time-consuming task for even an experienced design engineer. He or she will have to make a selection trial many times before ideal sufficient accuracy can be reached. The deficiency of the FEM becomes much explicit in the accurate computation.

The problems stem from the theoretical drawbacks of the FEM, which embody two major aspects. One is as follows, a single element defined domain of a node shape function contains only a portion of a full node shown in Fig.1, that means the construction process of a shape function in the traditional FEM can be viewed as a procedure in which a full node is broken up into split nodes, resulting in splitting all full nodes over a structural domain in meshing step and discretizing a whole structural domain. Another is that a mesh generation is usually an arbitrary operation based on the empiricism without a solid mathematical foundation to support, which leads to two optional schemes of regular or irregular meshing and could bring about a random node distribution over a structural domain. Thus, it is called an irrational MRA when the random node layout is allowed. The theoretical deficiencies of the standard FEM come from the element split node and the irrational MRA.

The great efforts have been made over the past thirty years to tackle the pitfalls of the FEM. The solutions were found and tracked in two trends. In order to make the meshing more efficient and accessible, the various improved FEM have been proposed, such as the automatic adaptive FEM [2, 3], the multigrid FEM [4, 5], etc. However, the source of the problems still remains. For eliminating the grid, investigators have come up with many fresh methods, such as the wavelet finite element method (WFEM)[6, 7], the meshfree method (MFM)[8, 9], the natural element method (NEM)[10, 11] and the isogeometric analysis method (IGAM) [12, 13],etc, which are featured with the full node shape function and the integrated structural computational model. Compared with the conventional FEM, these meshless methods have an advantage not to mesh over a structural domain. Therefore, these methods have illustrated their powerful capability and computational efficiency in dealing with some problems. However, they always have such major inherent deficiencies as the complexity of full node shape function construction by tensor-product or polynomial coefficient numerical simulating technique, the absence of the Kronecker delta property and the lack of a rigorous mathematical basis for node distribution over a structural domain (an irregular node configuration makes an irrational MRA), which make the treatment of boundary condition complicated and the selection of node layout empirical, that substantially reduce computational efficiency. Hence, these meshless methods have never found a wide application in engineering practice just as the FEM.

The drawbacks of all those above-mentioned methods can be eliminated by the introduction of a new multiresolution element method in this paper. With respect to the plate element in the finite element stock, a new multiresolution triangular plate-bending element is formulated by a new MRA, which is constituted by scaled and translated version as subspace basis functions of a basic full node shape function. The basic full node shape function is then constructed from making a series of parallelograms to superimpose identical triangle-defined domains around the origin of coordinates. As we can see, the full node shape function is quite simple, clear and has the Kronecker delta property. The mesh over a structural domain can be deleted and the node layout of uniform distribution pattern (the regular meshing scheme) can be reasonably defined because the proposed method possesses the rational MRA based on a simple, clear and rigorous mathematical basis, which endows the proposed element with the resolution level (RL) that can be modulated to freely change the node number and position at a request, adjusting structural analysis clarity accordingly. As a result, the selection of an appropriate RL is more timely efficient than the selection of an appropriate meshing scheme and the proposed element method can bring about substantial improvement of the computational efficiency in the structural analysis when compared with the corresponding FEM or other meshless methods.

## 2. Basic full node shape function

As shown in Fig.1., an arbitrary triangle plate element is set against a Carstesian coordinate system with the geometric configuration of the bottom sideline length as $a$, the height as $h$. Obviously, The analytical functions for the bottom sideline in the coordinate system can be written in dimensionless quantity as

$$\frac{y}{h}=0 \qquad (1)$$

For the other sideline (not one that goes through the coordinate origin) is assumed as

$$\frac{x}{a}+\frac{y}{b}=1 \qquad (2)$$

Where $a$,$b$ are denoted as the horizontal and the vertical intercepts respectively.
For the third is determined as

$$\frac{x}{a}-\left(\frac{1}{h}-\frac{1}{b}\right)y=0 \quad (b \geq h) \qquad (3)$$

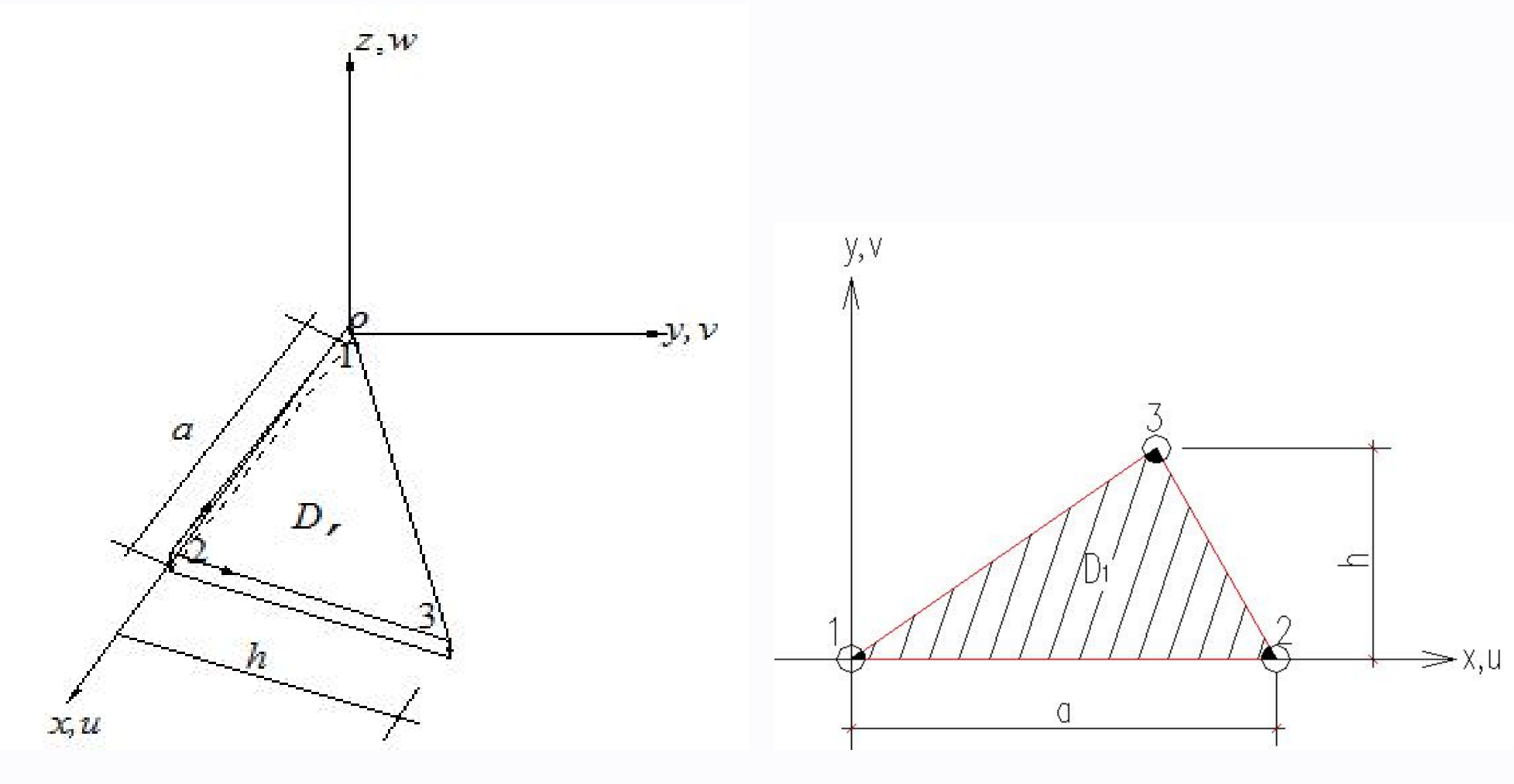

***a***. A triangular plate　　　　　　　***b*** Mid-plane of the plate

**Fig 1.** A triangular plate-bending element

Afterward, the transverse displacement $w^e$ in $z$ axis direction at an arbitrary point within the triangular plate element can be defined as

$$w^e=\sum_{i=1}^{3}N_i w_i+\sum_{i=1}^{3}N_{xi}\theta_{xi}+\sum_{i=1}^{3}N_{yi}\theta_{yi} \qquad (4)$$

where $w_i$ , $\theta_{xi}$ , $\theta_{yi}$ are the transverse, rotational displacements at node $i$ of the element respectively in the Carstesian coordinate system. $N_i$ , $N_{xi}$ , $N_{yi}$ are the conventional shape functions at the node $i$ ( $i$ =1,2,3), which are defined on the domain $D_1$ as follows

$$N_1 = L_1 + L_1^2 L_2 + L_1^2 L_3 - L_1 L_2^2 - L_1 L_3^2$$

$$N_{x1} = -b_3 \left( L_1^2 L_2 + \frac{1}{2} L_1 L_2 L_3 \right) + b_2 \left( L_3 L_1^2 + \frac{1}{2} L_1 L_2 L_3 \right)$$

$$N_{y1} = -c_3 \left( L_1^2 L_2 + \frac{1}{2} L_1 L_2 L_3 \right) + c_2 \left( L_3 L_1^2 + \frac{1}{2} L_1 L_2 L_3 \right)$$

$$b_3 = y_1 - y_2, \quad c_3 = x_2 - x_1, \quad b_2 = y_3 - y_1, \quad c_2 = x_1 - x_3$$

Based on the analytical functions for the three triangle sidelines obtained above, the following relationship can be gotten:

$$L_1 = 1 - \left( \frac{x}{a} + \frac{y}{b} \right), \quad L_2 = \frac{x}{a} - \left( \frac{1}{h} - \frac{1}{b} \right) y, \quad L_3 = \frac{y}{h}, \quad x, y \in D_1 \qquad (5a,b,c)$$

Obviously, there exists relationship $L_1 + L_2 + L_3 = 1$

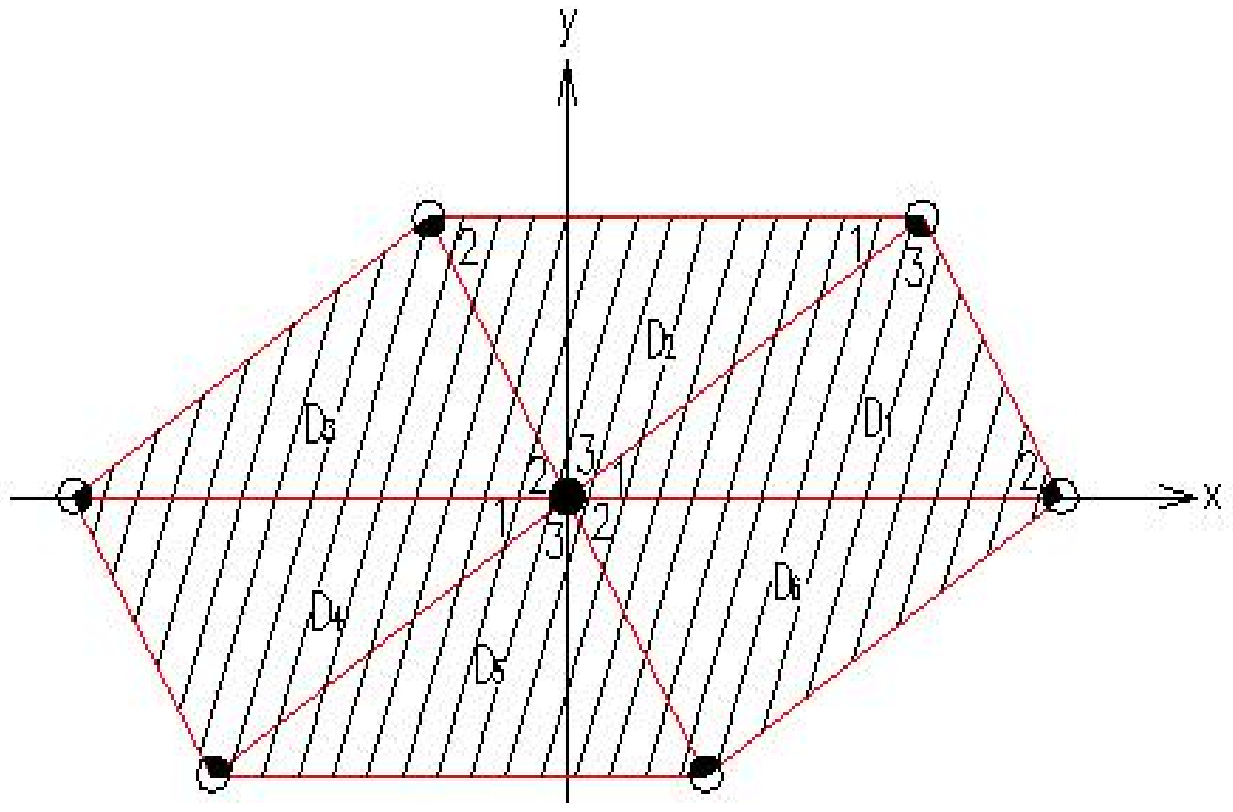

**Fig 2.** The extended hexagon domain enclosing a node at the coordinate origin

As we can see in Fig.1, the supporting domain (shaded area $D_1$) of the triangular element contains only a part (blackened portion) of a full node, that means a full node is broken up into split nodes in the process of the traditional node shape function construction and all full nodes within the structural domain are thus discretized by meshing. In order to formulate a full node shape function or constitute an integrated computational model, the single triangle-defined domain for the split node should be extended to the hexagon area by means of successively building up a series of parallelograms to superimpose identical triangle-defined regimes around the coordinate origin. Subsequently, the node at the coordinate zero is enclosed by the hexagon domain (shaded area) as shown in Fig.2. The basic shape function for the full node (blackened node) at the coordinate origin can be defined as following:

$$\phi(x,y)=\begin{cases} N_1(x,y) & (x,y\in D_1) \\ N_3(x,y) & (x,y\in D_2) \\ N_2(x,y) & (x,y\in D_3) \\ N_1(x,y) & (x,y\in D_4) \\ N_3(x,y) & (x,y\in D_5) \\ N_2(x,y) & (x,y\in D_6) \end{cases} \tag{6}$$

$$\phi_x(x,y)=\begin{cases} N_{x1}(x,y) & (x,y\in D_1) \\ N_{x3}(x,y) & (x,y\in D_2) \\ N_{x2}(x,y) & (x,y\in D_3) \\ N_{x1}(x,y) & (x,y\in D_4) \\ N_{x3}(x,y) & (x,y\in D_5) \\ N_{x2}(x,y) & (x,y\in D_6) \end{cases} \tag{7}$$

$$\phi_y(x,y)=\begin{cases} N_{y1}(x,y) & (x,y\in D_1) \\ N_{y3}(x,y) & (x,y\in D_2) \\ N_{y2}(x,y) & (x,y\in D_3) \\ N_{y1}(x,y) & (x,y\in D_4) \\ N_{y3}(x,y) & (x,y\in D_5) \\ N_{y2}(x,y) & (x,y\in D_6) \end{cases} \tag{8}$$

where $N_i$ , $N_{xi}$ , $N_{yi}$ are the shape functions at the node $i$ ( $i$ =1,2,3), which are defined on the domains of $D_1$, $D_2$ , $D_3$ , $D_4$ , $D_5$ , $D_6$ corresponding to six split nodes around the coordinate origin respectively.

In light of the regular node shape function construction method by area coordinates for a triangular plate element, the six split node shape functions can be founded by the analytical functions for the six sidelines of the hexagon in the Carstesian coordinate system. Based on the analytical functions for three triangular sidelines (1), (2), (3), the three upper hexagon sideline functions are easily written as

$$\frac{x}{a}+\frac{y}{b}=1 \tag{9a}$$

$$\left(\frac{1}{h}-\frac{1}{b}\right)y-\frac{x}{a}=1 \tag{9b}$$

$$\frac{y}{h}=1\ , \tag{9c}$$

Therefore, the three lower hexagon sideline function expressions can be easily obtained by shifting each upper sideline a distance along $x$,$y$ axis respectively, that is

$$\frac{x}{a}+\frac{y}{b}=-1 \tag{10a}$$

$$\left(\frac{1}{h}-\frac{1}{b}\right)y-\frac{x}{a}=-1 \tag{10b}$$

$$\frac{y}{h}=-1\ , \tag{10c}$$

As a result, based on the hexagon sideline analytical functions, the split node shape functions on the domains of $D_2$, $D_3$, $D_4$, $D_5$, $D_6$ can be founded respectively as

$$\begin{cases} L_1=\frac{x}{a}+\frac{y}{b}, \quad L_2=\left(\frac{1}{h}-\frac{1}{b}\right)y-\frac{x}{a}, \quad L_3=1-\frac{y}{h}, \quad x,y\in D_2 \\ N_3=L_3+L_3^2L_1+L_3^2L_2-L_3L_1^2-L_3L_2^2 \\ N_{x3}=-b_2\left(L_3^2L_1+\frac{1}{2}L_1L_2L_3\right)+b_1\left(L_2L_3^2+\frac{1}{2}L_1L_2L_3\right) \\ N_{y3}=-c_2\left(L_3^2L_1+\frac{1}{2}L_1L_2L_3\right)+c_1\left(L_2L_3^2+\frac{1}{2}L_1L_2L_3\right) \\ b_2=y_3-y_1, \quad c_2=x_1-x_2, \quad b_1=y_2-y_3, \quad c_1=x_3-x_2 \end{cases}$$

$$\begin{cases} L_1=-\left(\frac{x}{a}+\frac{y}{b}\right), \quad L_2=1-\left(\frac{1}{h}-\frac{1}{b}\right)y+\frac{x}{a}, \quad L_3=\frac{y}{h}, \quad x,y\in D_3 \\ N_2=L_2+L_2^2L_3+L_2^2L_1-L_2L_3^2-L_2L_1^2 \\ N_{x2}=-b_1\left(L_2^2L_3+\frac{1}{2}L_1L_2L_3\right)+b_3\left(L_1L_2^2+\frac{1}{2}L_1L_2L_3\right) \\ N_{y2}=-c_1\left(L_2^2L_3+\frac{1}{2}L_1L_2L_3\right)+c_3\left(L_1L_2^2+\frac{1}{2}L_1L_2L_3\right) \\ b_1=y_2-y_3, \quad c_1=x_3-x_2, \quad b_3=y_1-y_2, \quad c_3=x_2-x_1 \end{cases}$$

$$\begin{cases} L_1=1+\frac{x}{a}+\frac{y}{b}, \quad L_2=\left(\frac{1}{h}-\frac{1}{b}\right)y-\frac{x}{a}, \quad L_3=-\frac{y}{h}, \quad x,y\in D_4 \\ N_1=L_1+L_1^2L_2+L_1^2L_3-L_1L_2^2-L_1L_3^2 \\ N_{x1}=-b_3\left(L_1^2L_2+\frac{1}{2}L_1L_2L_3\right)+b_2\left(L_3L_1^2+\frac{1}{2}L_1L_2L_3\right) \\ N_{y1}=-c_3\left(L_1^2L_2+\frac{1}{2}L_1L_2L_3\right)+c_2\left(L_3L_1^2+\frac{1}{2}L_1L_2L_3\right) \\ b_3=y_1-y_2, \quad c_3=x_2-x_1, \quad b_2=y_3-y_1, \quad c_2=x_1-x_3 \end{cases}$$

$$\begin{cases} L_1 = -\left(\dfrac{x}{a}+\dfrac{y}{b}\right), \quad L_2 = \dfrac{x}{a}-\left(\dfrac{1}{h}-\dfrac{1}{b}\right)y, \quad L_3 = 1+\dfrac{y}{h}, \ x, y \in D_5 \\ N_3 = L_3 + L_3^2 L_1 + L_3^2 L_2 - L_3 L_1^2 - L_3 L_2^2 \\ N_{x3} = -b_2\left(L_3^2 L_1 + \dfrac{1}{2}L_1 L_2 L_3\right) + b_1\left(L_2 L_3^2 + \dfrac{1}{2}L_1 L_2 L_3\right) \\ N_{y3} = -c_2\left(L_3^2 L_1 + \dfrac{1}{2}L_1 L_2 L_3\right) + c_1\left(L_2 L_3^2 + \dfrac{1}{2}L_1 L_2 L_3\right) \\ b_2 = y_3 - y_1, \quad c_2 = x_1 - x_2, \quad b_1 = y_2 - y_3, \quad c_1 = x_3 - x_2 \end{cases}$$

$$\begin{cases} L_1 = \dfrac{x}{a}+\dfrac{y}{b}, \quad L_2 = 1-\dfrac{x}{a}+\left(\dfrac{1}{h}-\dfrac{1}{b}\right)y, \quad L_3 = -\dfrac{y}{h}, \ x, y \in D_6 \\ N_2 = L_2 + L_2^2 L_3 + L_2^2 L_1 - L_2 L_3^2 - L_2 L_1^2 \\ N_{x2} = -b_1\left(L_2^2 L_3 + \dfrac{1}{2}L_1 L_2 L_3\right) + b_3\left(L_1 L_2^2 + \dfrac{1}{2}L_1 L_2 L_3\right) \\ N_{y2} = -c_1\left(L_2^2 L_3 + \dfrac{1}{2}L_1 L_2 L_3\right) + c_3\left(L_1 L_2^2 + \dfrac{1}{2}L_1 L_2 L_3\right) \\ b_1 = y_2 - y_3, \quad c_1 = x_3 - x_2, \quad b_3 = y_1 - y_2, \quad c_3 = x_2 - x_1 \end{cases}$$

in which there all exist relationships

$1 \ge L_1 \ge 0, \quad 1 \ge L_2 \ge 0, \quad 1 \ge L_3 \ge 0, L_1 + L_2 + L_3 = 1$ for the various domains of $D_1$, $D_2$, $D_3$, $D_4$, $D_5$, $D_6$ respectively.

Up to now, the basic full node shape functions $\phi(x,y)$, $\phi_x(x,y)$, $\phi_y(x,y)$ can be graphed respectively in Fig.3.

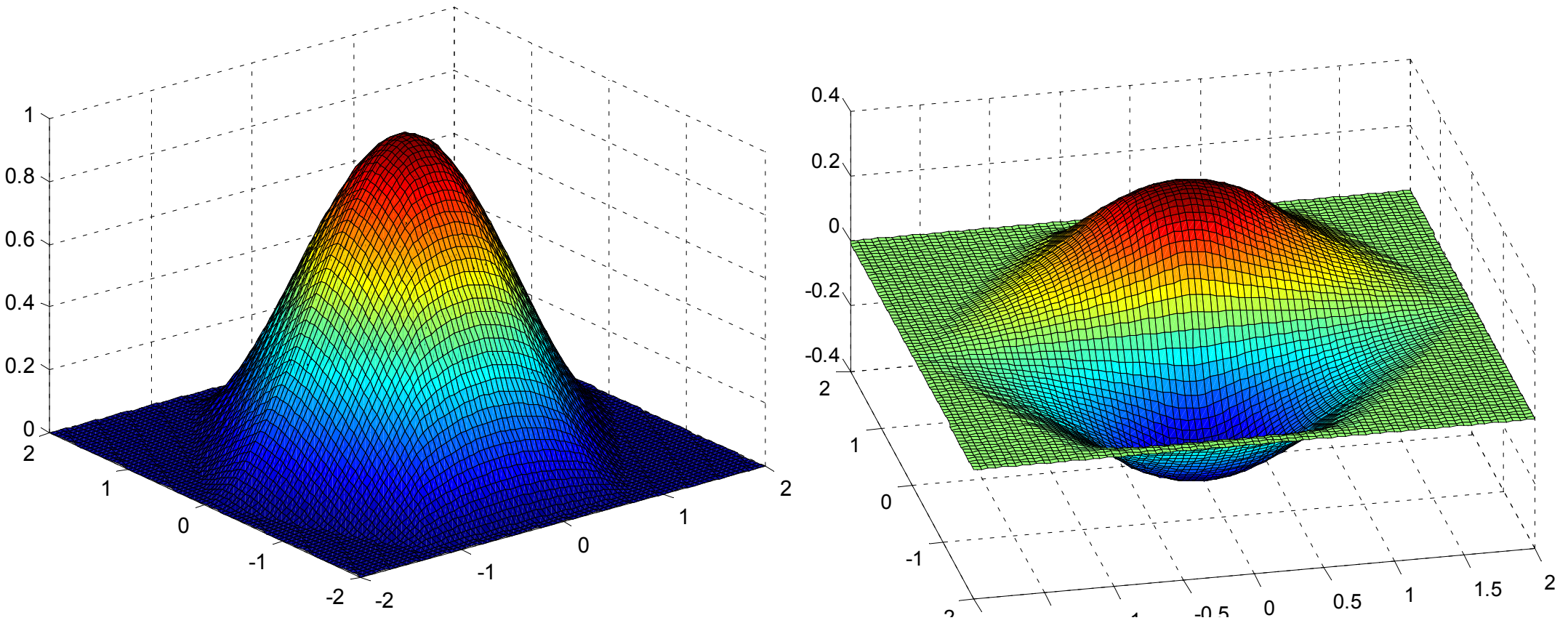

*a.* $\phi(x,y)$ *b.* $\phi_x(x,y)$

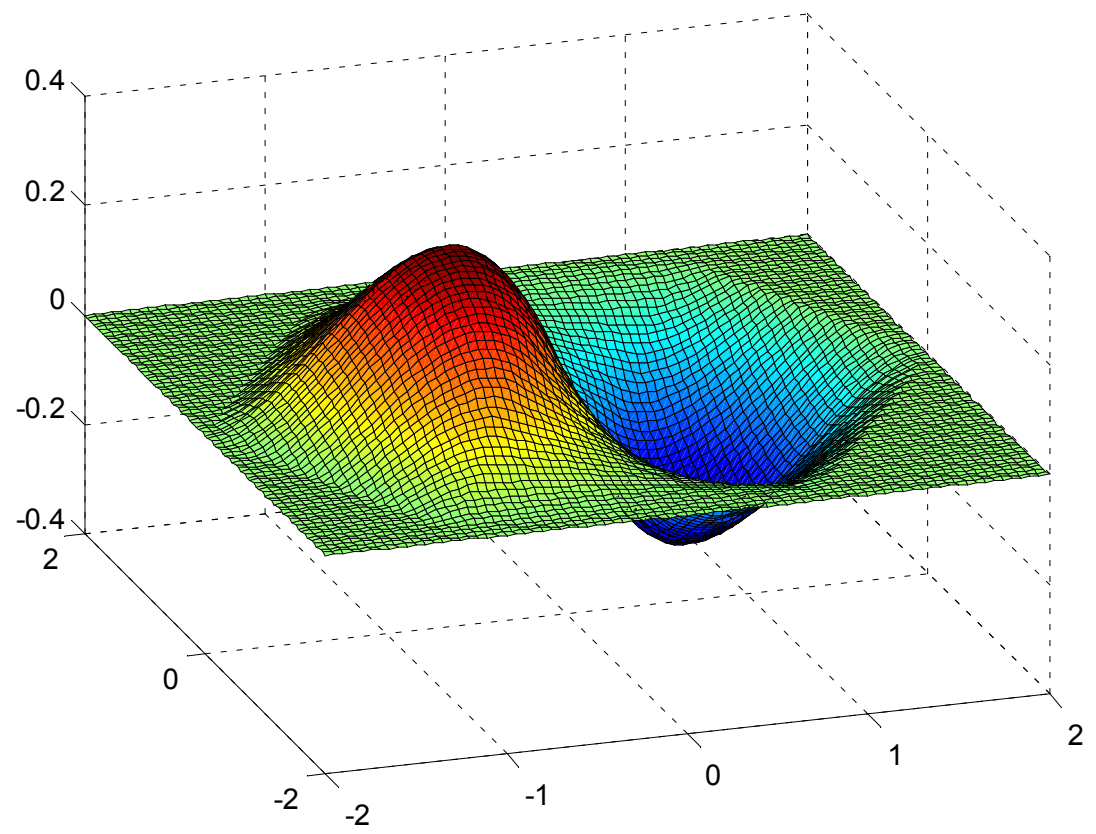

c. $\phi_y(x,y)$

**Fig.3.**The basic full node shape functions $\phi(x,y)$, $\phi_x(x,y)$, $\phi_y(x,y)$ at the coordinate origin

It is evident that the basic full node shape functions $\phi(x,y), \phi_x(x,y), \phi_y(x,y)$ are continuous and possess the Kronecker delta property as follows:

$$\begin{cases} \phi(0,0)=1 \quad \phi(6\,points)=0 \quad \dfrac{\partial\phi(0,0)}{\partial x}=0 \\ \dfrac{\partial\phi(0,0)}{\partial y}=0 \quad \dfrac{\partial\phi(6\,points)}{\partial x}=0 \quad \dfrac{\partial\phi(6\,points)}{\partial y}=0 \\ \phi_x(0,0)=0 \quad \phi_x(6\,points)=0 \quad \dfrac{\partial\phi_x(0,0)}{\partial x}=0 \\ \dfrac{\partial\phi_x(0,0)}{\partial y}=1 \quad \dfrac{\partial\phi_x(6\,points)}{\partial x}=0 \quad \dfrac{\partial\phi_x(6\,points)}{\partial y}=0 \\ \phi_y(0,0)=0 \quad \phi_y(6\,points)=0 \quad \dfrac{\partial\phi_y(0,0)}{\partial x}=-1 \\ \dfrac{\partial\phi_y(0,0)}{\partial y}=0 \quad \dfrac{\partial\phi_y(6\,points)}{\partial x}=0 \quad \dfrac{\partial\phi_y(6\,points)}{\partial y}=0 \end{cases} \tag{11}$$

## 3. Displacement Subspace Sequence

In order to carry out a MRA of a thin plate structure, the mutual nesting displacement subspace sequence for a plate element should be established. In this paper, a totally new technique is proposed to construct the MRA which is based on the concept that a subspace sequence (multi-resolution subspaces) can be formulated by subspace basis function vectors at different resolution levels whose elements-scaling function vector can be constructed by scaling and shifting on the domain of the basic full node shape functions. As a result, the displacement subspace basis function vector at an arbitrary resolution level (RL) of $\frac{1}{2}(m+1)\times(m+2)$ for a triangular plate-bending element is formulated as follows:

$$\mathbf{\Psi}_m = \left[ \mathbf{\Phi}_{mm,00} \cdots \mathbf{\Phi}_{mm,rs} \cdots \mathbf{\Phi}_{mm,mm} \right] \tag{12}$$

where $\mathbf{\Phi}_{mn,rs} = \left[ \phi_{mm,rs} \quad \phi_{xmm,rs} / m \quad \phi_{ymm,rs} / m \right]$ is the scaling basis function vector, $\phi_{mm,rs}=\phi(mx-\left(r-s\frac{h}{b}\right)a, my-sh)$ , $\phi_{xmm,rs}=\phi_x(mx-\left(r-s\frac{h}{b}\right)a, my-sh)$ , $\phi_{ymm,rs}=\phi_y(mx-\left(r-s\frac{h}{b}\right)a, my-sh)$ , $m$ denoted as the positive integers, the scaling parameters in $x$ ,$y$ directions respectively. $r$, $s$ as the positive integers, the node position parameters, that is $r=0,1,2,3\cdots m$ , $s=0,1,2,3\cdots m$ , Here, $m \ge r \ge s$ , $\left(mx-ra\right) \in \left[-a,a\right]$ , $\left(my-sh\right) \in \left[-h,h\right]$ , $x, y \in D_1$ .

It is seen from Eq. (12) that the nodes for the scaling process are equally spaced on the triangle domain $D_1$ in $x$ ,$y$ directions respectively

Scaling of the basic full node shape function on the domain of $\left[-a,a\right]\times\left[-h,h\right]$ ( precisely on the domain of $\left[-\frac{a}{m},\frac{a}{m}\right]\times\left[-\frac{h}{m},\frac{h}{m}\right]$) and then shifting to other nodes $\left(\frac{a}{m}\left(r-s\frac{h}{b}\right),\frac{s}{m}h\right)$ within the element domain $D_1$ $\left(m \ge r \ge s\right)$ will produce the various full node shape functions..

As shown in the fig.4., for a full node within the triangular element domain, a scaled and shifted version of the basic full node shape function is adopted; For a 1/2 full node denoted as B on the triangular element boundary, a half of the scaled and shifted version is employed (the other half out of the triangular element domain is deleted); For a split node denoted as T on the tip, only 1 split node shape function in the scaled and shifted version is applied (the other 5 splits out of the triangular element domain are cancelled).

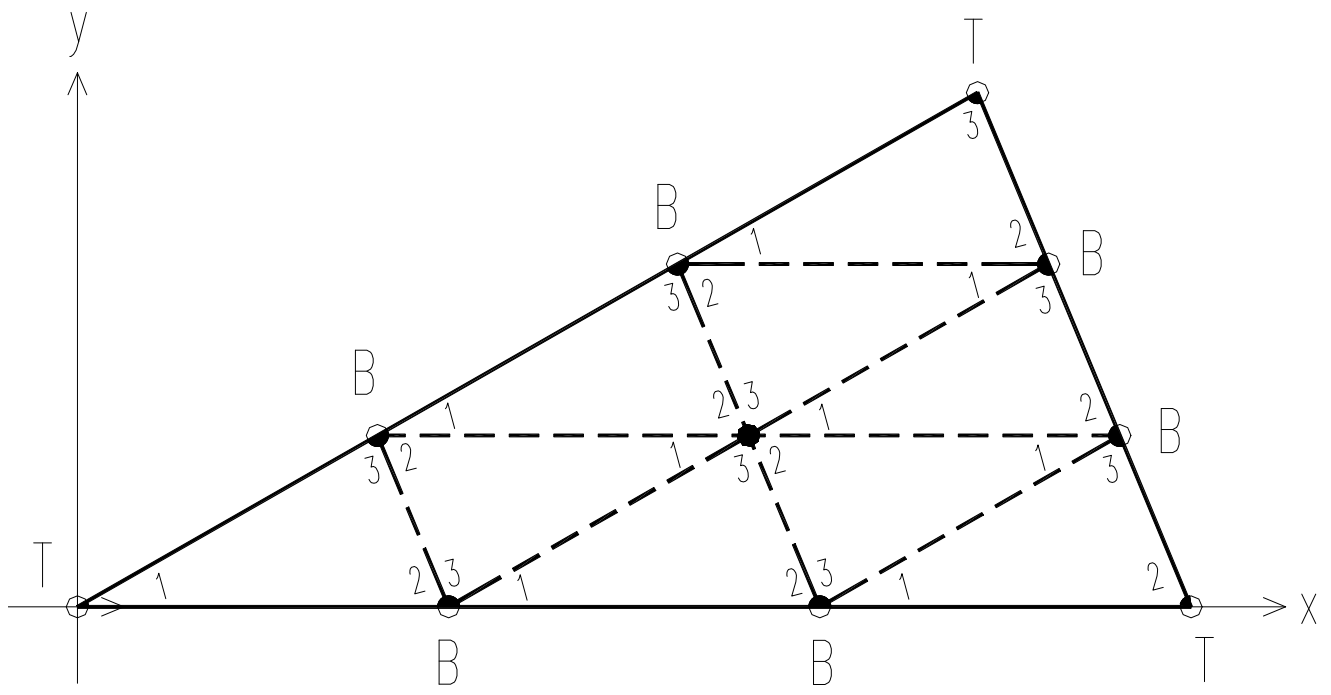

**Fig.4.**The node partition for a triangular element when the RL=2×5 (m=3)

Since the elements in the base functions are linearly independent with the various scaling and the different shifting parameters, the subspaces in the subspace sequence can be established, thus formulating a MRA framework, that is

$$\mathbf{W}_m = \left[V_1 \cdots \quad V_i \cdots \quad V_m\right] \tag{13a}$$

$$V_i := span\left\{\mathbf{\Psi}_i : i \in Z\right\} \tag{13b}$$

If $I = 2i$, then $V_i \subset V_I$ (13c)

Thus, it can be found that the displacement subspace sequence $\mathbf{W}_m$ can be taken for a solid mathematical foundation for the MRA framework and $V_1$ is equivalent to the displacement field for a traditional 3-node plate-bending triangular element that is the reason why the traditional triangular plate element is regarded as a mono-resolution one and also a special case of the multiresolution triangular.

Based the MRA established, the deflection of the triangular plate element in the displacement subspace at RL of $\frac{1}{2}(m+1)\times(m+2)$ can be defined as follows

$$w_m^e = \mathbf{\Psi}_m \mathbf{a}_m^e \tag{14}$$

where $\mathbf{a}_{mn}^e = \left[\left[w_{00}, \theta_{x00}, \theta_{y00}\right] \cdots \quad \left[w_{rs}, \theta_{xrs}, \theta_{yrs}\right] \cdots \quad \left[w_{mn}, \theta_{xmn}, \theta_{ymn}\right]\right]^T$, $w_{rs}, \theta_{xrs}, \theta_{yrs}$ are the transverse and rotational displacements respectively at the element node $\left(\frac{a}{m}\left(r - s\frac{h}{b}\right), \frac{s}{m}h\right)$.

It is obvious that the proposed multi-resolution element is a meshfree one whose nodes are uniformly scattered, node number and position fully determined by the *RL*. When the scaling parameter m=1(*RL*=$\frac{1}{2}(2)\times(3) = 1\times 3$ ), that is a traditional 3-node triangular plate-bending element, eq. (14) will be reduced to eq. (4).

## 4 Multiresolution triangular plate element

The generalized function of potential energy in a displacement subspace at the resolution level of $\frac{1}{2}(m+1)\times(m+2)$ for a triangular plate-bending element can be defined as

$$\Pi(V_m) = \frac{1}{2}\iint_{D_1} [\kappa]_m^T [D_b][\kappa]_m \, dxdy - \iint_{D_1} q w_m^e dxdy - \sum_i Q_i w_{mi}^e \tag{15}$$

where $[\kappa]_m = -\begin{bmatrix} \dfrac{\partial^2 w_m^e}{\partial x^2} \\ \dfrac{\partial^2 w_m^e}{\partial y^2} \\ 2\dfrac{\partial^2 w_m^e}{\partial x \partial y} \end{bmatrix}$, $[D_b] = C_b \begin{bmatrix} 1 & \mu & 0 \\ \mu & 1 & 0 \\ 0 & 0 & (1-\mu)/2 \end{bmatrix}$, $C_b = \dfrac{Et^3}{12(1-\mu^2)}$, $E$ is the material Young modulus, $t$ the thickness of the element, $\mu$ the Poisson's ratio, $q$ distributed transverse loadings, $Q$ the lump transverse loadings elastic modulus.

$$[\kappa]_m = [B_{00}, \cdots B_{rs}, \cdots B_{mm}]\mathbf{a}_m^e \tag{16}$$

where

$$B_{rs} = -\begin{bmatrix} \dfrac{\partial^2 \phi_{rs,m}}{\partial x^2} & \dfrac{\partial^2 \phi_{xrs,m}}{\partial x^2} & \dfrac{\partial^2 \phi_{yrs,m}}{\partial x^2} \\ \dfrac{\partial^2 \phi_{rs,m}}{\partial y^2} & \dfrac{\partial^2 \phi_{xrs,m}}{\partial y^2} & \dfrac{\partial^2 \phi_{yrs,m}}{\partial y^2} \\ 2\dfrac{\partial^2 \phi_{rs,m}}{\partial x \partial y} & 2\dfrac{\partial^2 \phi_{xrs,m}}{\partial x \partial y} & 2\dfrac{\partial^2 \phi_{yrs,m}}{\partial x \partial y} \end{bmatrix}$$

Substituting (14) into (15) and consolidating, we get

$$\Pi_p(V_m) = \frac{1}{2}\mathbf{a}_m^{eT}\mathbf{K}_m^e\mathbf{a}_m^e - \mathbf{a}_m^{eT}\mathbf{f}_m^e - \mathbf{a}_m^{eT}\mathbf{F}_m^e \tag{17}$$

in which $\mathbf{K}_m^e$ is denoted as the element stiffness matrix; $\mathbf{f}_m^e$ as the element distributed equivalent node force vector; $\mathbf{F}_m^e$ as the element concentrated equivalent node force vector.

According to the principal of minimums potential energy $\delta\Pi_p(V_m) = 0$, the following element equilibrium equations can be obtained

$$\mathbf{K}_m^e\mathbf{a}_m^e = \mathbf{f}_m^e + \mathbf{F}_m^e$$

where $\mathbf{K}_m^e = \begin{bmatrix} \mathbf{k}_{00}^{00} \cdots & \mathbf{k}_{rs}^{00} \cdots & \mathbf{k}_{mm}^{00} \\ \cdot & \cdot & \cdot \\ \cdot & \cdot & \cdot \\ \cdot & \cdot & \cdot \\ \mathbf{k}_{00}^{rs} \cdots & \mathbf{k}_{rs}^{rs} \cdots & \mathbf{k}_{ij}^{rs} \\ \cdot & \cdot & \cdot \\ \cdot & \cdot & \cdot \\ \cdot & \cdot & \cdot \\ \mathbf{k}_{00}^{mm} \cdot & \mathbf{k}_{rs}^{mm} \cdot & \mathbf{k}_{mm}^{mm} \end{bmatrix}$

in which the superscript denoted as the row number of the matrix and the subscript as the aligned element node numbering (*r*, *s*). In terms of the properties of the extended shape functions, we have

$$\begin{cases} \mathbf{k}_{rs}^{rs} = \sum\limits_{\substack{|c-r|\le 1 \\ |d-s|\le 1}} \mathbf{k}_{cd,rs} \\ \mathbf{k}_{rs}^{rs} = \mathbf{k}_{cd,rs} = 0, when\, |c-r| > 1, |d-s| > 1 \end{cases} \tag{18}$$

in which $\mathbf{k}_{cd,rs}$ is the coupled node stiffness matrix relating the node (*c*, *d*) to (*r*, *s*).

$$\mathbf{k}_{cd,rs} = \iint\limits_{D_1} \mathbf{B}_{cd} \mathbf{E} \mathbf{B}_{rs}^{b} dxdy \tag{19}$$

$$\begin{cases} \mathbf{f}_{m,rs}^{er} = \iint \left[ \mathbf{\Phi}_{mm,rs} \right]^T q dxdy \\ \mathbf{F}_{m,rs}^{er} = \sum\limits_{i} \left[ \mathbf{\Phi}_{mm,rs} \left( mx_i - ra, my_i - sh \right) \right]^T P_i \end{cases} \tag{20}$$

where $x_i$, $y_i$ is the local coordinate at the locations the lump loading acting on.

## 5.Transformation Matrix

In order to carry out structural analysis, the element stiffness and mass matrices $\mathbf{K}_m^e$ the loading column vectors $\mathbf{f}_m^e$ , $\mathbf{F}_m^e$ should be transformed from the element local coordinate system (*xoy*) to the structural global coordinate system (*XOY*). The transforming relations from the local to the global are defined as follows:

$$\mathbf{K}_m^i = \mathbf{T}_m^{eT} \mathbf{K}_m^e \mathbf{T}_m^e \tag{21}$$

$$\mathbf{f}_m^i = \mathbf{T}_m^{eT} \mathbf{f}_m^e \tag{22}$$

$$\mathbf{F}_m^i = \mathbf{T}_m^{eT} \mathbf{F}_m^e \tag{23}$$

where $\mathbf{K}_m^i$ is the element stiffness matrix, $\mathbf{f}_m^i$ , $\mathbf{F}_m^i$ the element loading column vectors under the global coordinate system. $\mathbf{T}_m^e$ is the element transformation matrix defined as follows;

$$\mathbf{T}_m^e = \begin{bmatrix} \boldsymbol{\lambda}_{11} & & & & \mathbf{0} \\ & \cdots & & & \\ & & \boldsymbol{\lambda}_{ij} & & \\ & & & \cdots & \\ \mathbf{0} & & & & \boldsymbol{\lambda}_{mm} \end{bmatrix} \quad \boldsymbol{\lambda}_{ij} = \begin{bmatrix} \cos\theta_{zZ} & 0 & 0 \\ 0 & \cos\theta_{xX} & \cos\theta_{xY} \\ 0 & \cos\theta_{yX} & \cos\theta_{yY} \end{bmatrix}$$

where $\theta$ is the intersection angle between the local and the global coordinate axes . The structural global stiffness, mass matrix $\mathbf{K}_m$, the global loading column vectors $\mathbf{f}_m$, $\mathbf{F}_m$ can be obtained by splicing

# 6. Numerical Example and Discussion

## 6.1 Numerical Example

**Example 1.** A simply supported square plate with the geometric configuration of length *L*, the thickness *t*, the Poisson's ratio $\mu$ and subjected to the uniform loading *q*. Evaluate the deflection and the bending moment at the center point of the plate.

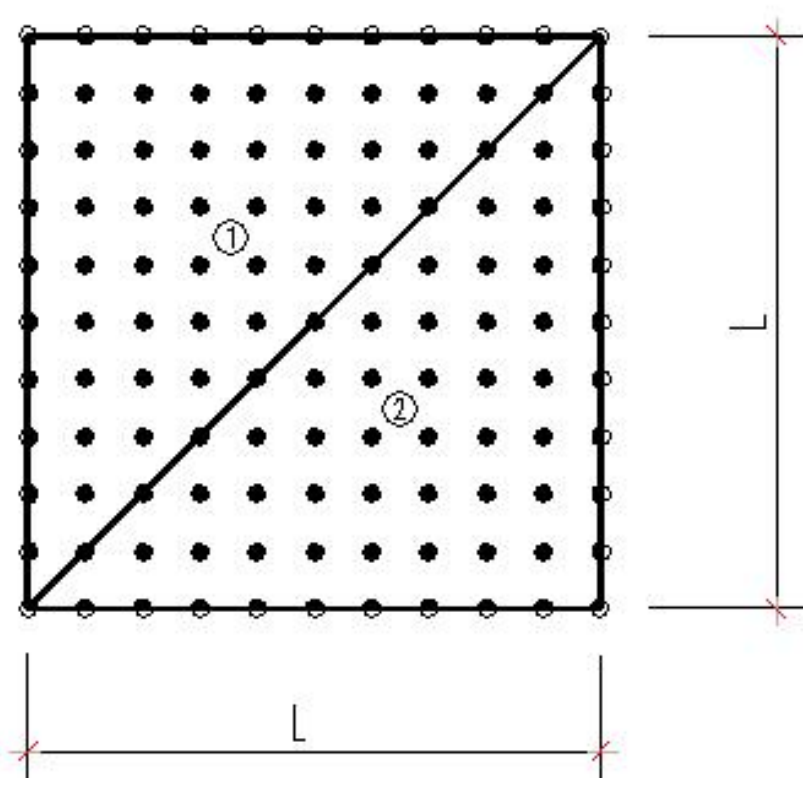

(a) The multiresolution element model (an integrated model by RL)

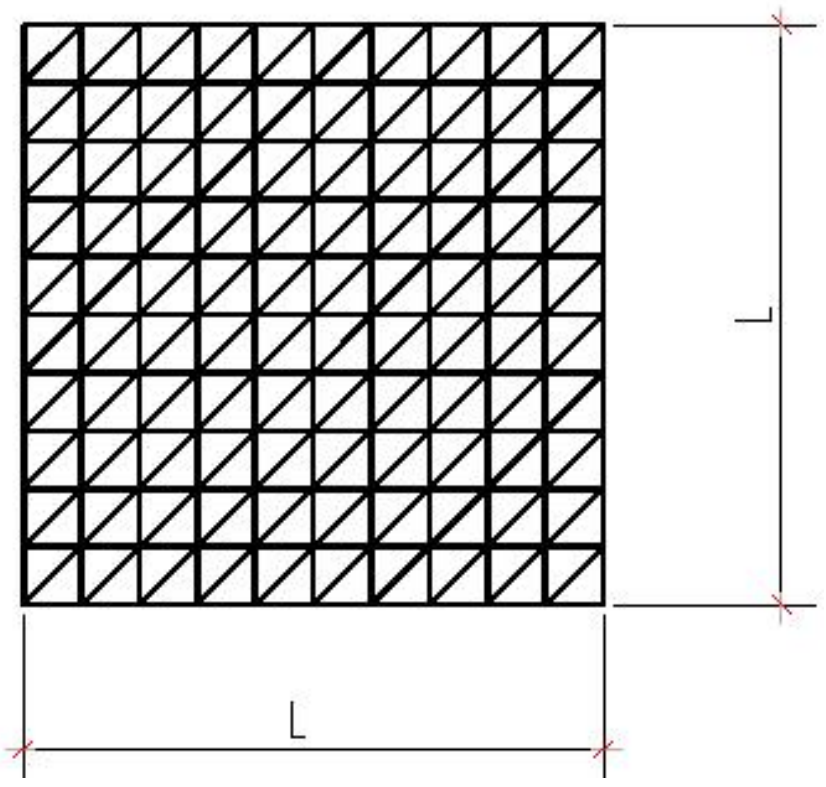

(b) The monoresolution element model (a discretized model by mesh)

**Fig.5**. Finite element models for the square plate

The problem is usually tackled with regular meshes typified by that shown in Fig.5(b) These meshes are built of two-triangle rectangular mesh units and identified as $N_x \times N_y$ , which denote the number of subdivisions in the *x* and the *y* dimensions, It is apparent that the meshing scheme is optional (irregular meshing scheme). The solution can also be found by the multiresolution triangular plate elements typified by that shown in Fig.5(a). These nodes are uniformly scattered in the two multiresolution triangular elements ①,② denoted by the RL, which defines a uniform node distribution pattern (regular meshing scheme). In the analysis process, these two multiresolution elements are spliced together along the common intersection.

The displacement responses and moments are evaluated by the proposed multiresolution element method, the traditional monoresolution element method. The central deflections and the bending moments for the plate with the boundary conditions as simply supported four edges (SS) and the clamped four edges (SC) under the different RLs and meshes are displayed in Table 1. The RL of the proposed and the corresponding meshes of the conventional are compared. It can be seen that the analysis accuracies with the proposed and the conventional are gradually improved respectively with the RL reaching high and the mesh approaching dense. However, the RL adjusting is more rationally and efficiently to be implemented than the meshing and remeshing to modulate element node number because the RL adjusting is based on the MRA framework which is constructed on a rigorous mathematical basis while the meshing or remeshing, which resorts to the empiricism, has no MRA framework. Thus, the computational efficiency of the proposed element method is higher than the traditional one. In this way, the proposed element exhibits its strong capability of accuracy adjustment and its high power of resolution to identify details (nodes)of deformed structure by means of modulating its resolution level, just as a multi-resolution camera with a pixel in its taken photo as a node in the proposed element. There appears no mesh in the proposed element just as no grid in the image. Hence, an element of superior analysis clarity surely has more nodes when compared with that of the inferior just as a clearer photo contains more pixels. Due to the basic full node shape function, the stiffness matrix and the loading column vectors of a proposed element can be automatically acquired through quadraturing around nodes in the element matrix formation step while those of the traditional plate element obtained through complex artificially reassembling of the element matrix around the node-related elements in the re-meshing process for their split nodes in a conventional element, which contributes a lot to computation efficiency improvement of the proposed method. Moreover, since the multiresolution triangular plate element model of a structure usually contains much less elements than the traditional element model, thus requiring much less times of transformation matrix multiplying, the computation efficiency of the proposed method appears much higher than the traditional in the step of element matrix transformation. In addition, because of the simplicity and clarity of a full node shape function formulation with the Kronecker delta property and the solid mathematical basis for the new MRA framework, the proposed method is also superior to other corresponding MRA methods in terms of the computational efficiency, the application flexibility and extent. Hence, taking all those causes into account, the conclusion can be drawn that the multi-resolution triangular plate-bending element method is more rationally, easily and efficiently to be executed, when compared with the traditional triangular plate element method.

**Table.1**.the Central deflection and the Bending moment for the plate under different RLs and meshes

| RL/elem | Mesh | Central deflection($ql^4/100\ D_b$) | | | | Bending moment ($ql^2/10$) | | | |
|---|---|---|---|---|---|---|---|---|---|
| | | SS | | SC | | SS (central) | | SC(side middle) | |
| | | Multi | Mono | Multi | Mono | Multi | Mono | Multi | Mono |
| 2×3 | 2×2 | 0.3950 | 0.3950 | 0.0998 | 0.0998 | 0.5026 | 0.5026 | -0.3551 | -0.3551 |
| 3×5 | 4×4 | 0.4039 | 0.4039 | 0.1194 | 0.1194 | 0.4880 | 0.4880 | -0.4761 | -0.4761 |
| 5×9 | 8×8 | 0.4058 | 0.4058 | 0.1249 | 0.1249 | 0.4824 | 0.4824 | -0.5028 | -0.5028 |
| 9×17 | 16×16 | 0.4062 | 0.4062 | 0.1262 | 0.1262 | 0.4800 | 0.4800 | -0.5104 | -0.5104 |
| Analytical[14] | | 0.4062 | | 0.1265 | | 0.4789 | | -0.5133 | |

**Example2.** As shown in Fig.6, a two opposite edge simply supported and other two free $60^0$ skew plate with the geometric configuration of length $L$ and the Poisson's ratio $\mu = 0.3$ is subjected to the uniform transverse loading of magnitude $q$. Evaluate the deflection at the center point of the plate.

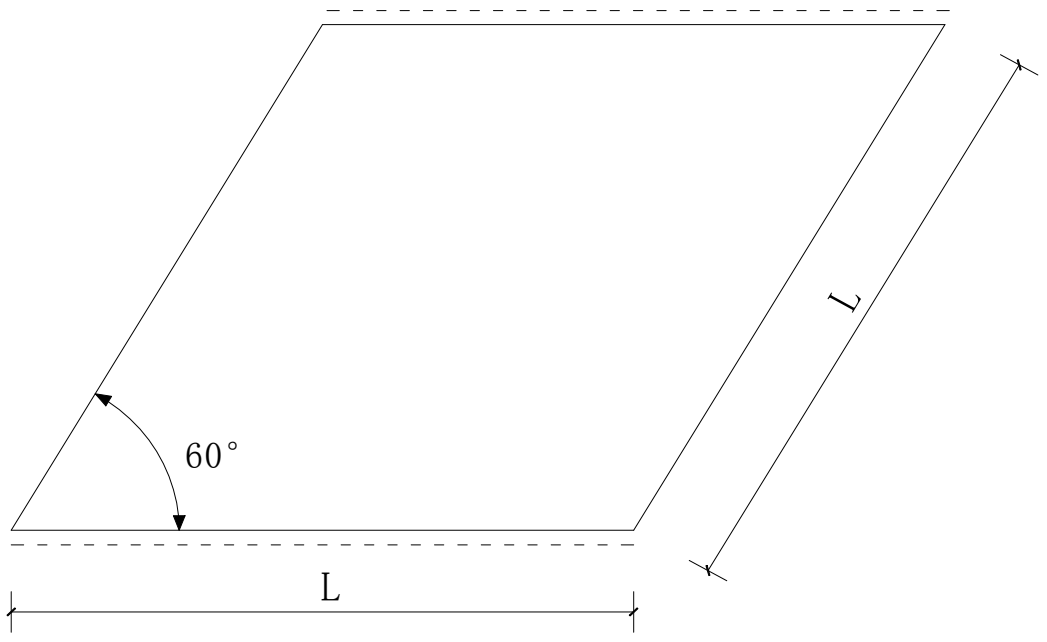

**Fig. 6.** A skew plate

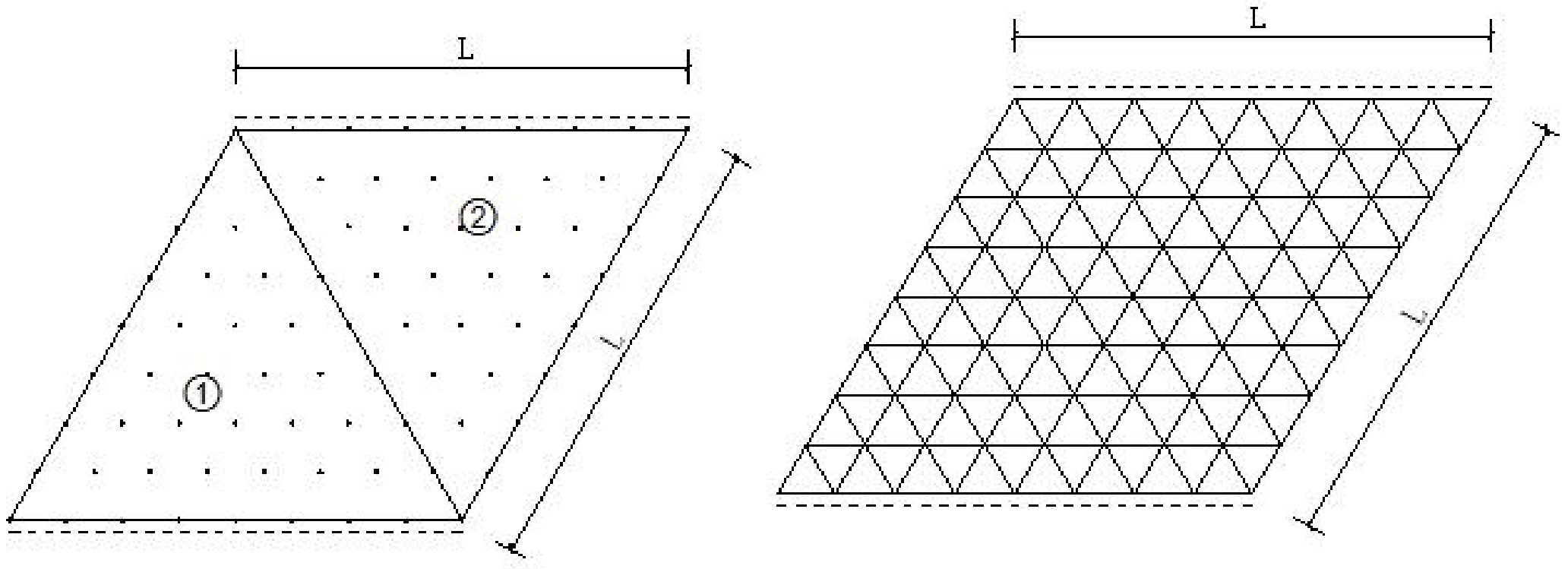

|(***a***) The multiresolution integrated model (***b***) The monoresolution discretized model

**Fig. 7** The finite element model for the skew plate

The problem is usually tackled with regular meshes typified by a discretized model (split node) shown in Fig.7(b) These meshes are built of two-triangle rectangular mesh units and identified as $N_x \times N_y$, which denote the number of subdivisions in the $x$ and the $y$ dimensions. The solution can also be found by the multiresolution triangular-bending elements typified by an integrated model (full node) shown in Fig.7(a). These nodes are uniformly scattered in the two multiresolution triangular plate-bending elements ①,② denoted by the RL. In the analysis process, these two multiresolution elements are spliced together along the common intersection boundary and the analysis clarity can be modulated by means of adjusting the RL. In addition, the wavelet model (full node) is made up of one 2D BSWI (B-spline wavelet on the interval) element of the jth scale=3, the mth order =4 are also employed abbreviated as BSWI43 with the DOF of $11 \times 11$. The results are summarized in the table2.It can be seen that the conventional and the proposed element methods exhibit identical monotonic increasing convergence to the 'exact' value with consistent mesh refinement and corresponding RL adjustment respectively.

**Table.2**. the center point deflection ( $\alpha \times qL^4/100D_0$ )

| Element type | | |
|---|---|---|
| The proposed RL/RL/elem | The conventional (mesh) | Deflection( α ) |
| 9×9/5×9 | 8×8 | 0.7920 |
| 13×13/7×13 | 12×12 | 0.7937 |
| 17×17/9×17 | 16×16 | 0.7930 |
| One BSWI [7] | 0.7925 | |
| Analytical [15] | 0.7945 | |

Although the BSWI43 is of high accuracy, when compared with the proposed, the deficiencies of the BSWI element are obvious as follows. In light of tensor product formulation of the multidimensional MRA framework, the DOF of a multi-dimensional BSWI element will be so drastically increased from that of a one-dimensional element in an irrational way, resulting in complex shape functions and substantial reduction of the computational efficiency. Secondly, due to the absence of Kronecker delta property of the tensor-product constructed shape functions, the special treatments should be taken to deal with the element boundary condition, which will bring about low computational efficiency. Thirdly, there exists no such a parameter as the RL with a clear mathematical sense. In addition, the RLs of the proposed and the corresponding meshes of the conventional are displayed in Table1. It can be found that the analysis clarities with the proposed and the conventional are gradually improved respectively with the RL reaching high and the mesh approaching dense. However, the RL adjusting is more rationally and efficiently to be implemented than the meshing and the re-meshing for the following two reasons. Firstly, the RL adjusting is based on the MRA framework that is constructed on a solid mathematical basis while the meshing or remeshing, which resorts to the empiricism, has no MRA framework. Secondly, the stiffness matrix and the loading column vectors of the proposed element can be obtained automatically around the nodes while those of the traditional triangular plate elements obtained by the artificially complex reassembling around the elements. Thus, the computational efficiency of the proposed element method is higher than the traditional one. In this way, the proposed plate element exhibits its strong capability of accuracy adjustment and its high power of resolution to identify details (nodes) of deformed structure by means of modulating its resolution level, just as a multiresolution camera with a pixel in its taken photo as a node in the proposed element. There appears no mesh in the proposed element just as no grid in a photo. Thus, an element of superior analysis accuracy surely has more nodes when compared with that of the inferior just as a clearer photo contains more pixels.

**Example 3.** A circular slab is subjected to the uniform transverse loading $q$ with its boundary conditions as: the edge is free or is fully clamped, and its geometry and physical parameters as: the radius $r$, the thickness $t$, the elasticity modulus $E$, the Poisson's ratio $\mu$=0.3. Find the displacement and bending moment at the central point of the slab.

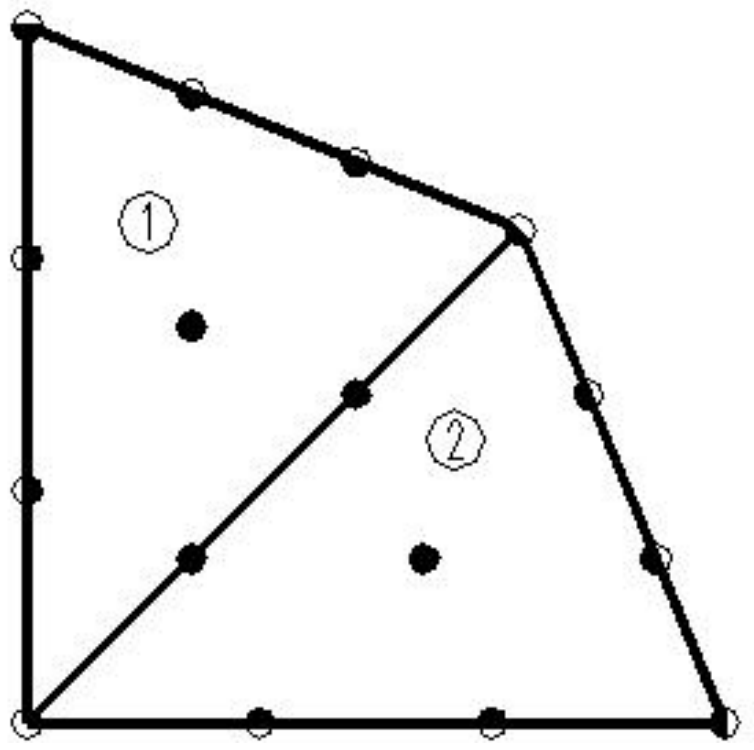

***a*** A.multiresolution integrated model

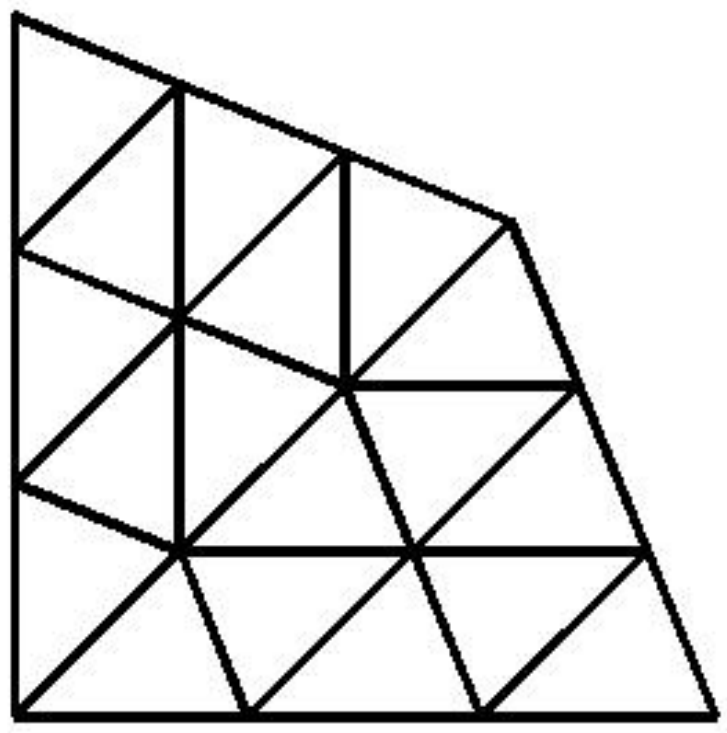

***b*** A monoresolution discretized model

**Fig.8.** The finite element model for the 1/4 circular slab

To calculate the displacement responses, symmetry conditions are exploited and The solution is done on the quadrant of the circular plate. the mutiresolution integrated model is herein composed of two multiresolution triangular plate-bending elements ①,②with each element RL of $2\times5$, hence the 1/4 slab RL of $4\times4$, as shown in Fig. 8*a*., and the monoresolution discretized one composed of the mesh of $3\times6$ displayed in Fig.8*b*. In the analysis process, these two multiresolution elements are spliced together along the common intersection boundary and the analysis clarity can be modulated by means of adjusting the RL. With respect to the conventional monoresolution, the structure is meshed into a group of monoresolution elements and the analysis accuracy is improved only by means of re-meshing. It can be seen that the RL adjusting is more rationally and easily to be implemented than the re-meshing because the proposed multiresolution element model of the circular plate structure contains much less elements than the monoresolution, hence requiring much less times of the transformation matrix multiplying, which results in much higher computational efficiency for the proposed element method than that for the traditional element method. The maximum displacement is summarized in Table.3.

**Table.3**.the Central deflection and Bending moment for the circular plate

| RL/elem | Mesh | Central deflection ( $\alpha \times ql^4/100 D_0$ ) | | | | Bending moment ( $\beta \times ql^2/10$) | |
|---|---|---|---|---|---|---|---|
| | | SS | | SC | | SC (central) | |
| | | Multi | Mono | Multi | Mono | Multi | Mono |
| 2×5 | 3×6 | 0.0145 | 0.0145 | 0.0638 | 0.0638 | 0.2103 | 0.2103 |
| 4×7 | 7×13 | 0.0153 | 0.0153 | 0.0637 | 0.0637 | 0.2073 | 0.2073 |
| Analytical[16] | | 0.0156 | | 0.0637 | | 0.20625 | |

## 6.2 Discussion

The RL can be referred to an ability to distinguish the magnitude difference in a distance between two points that is measured at a lengthy measurement, such as meter, centimeter or Nano scales etc. respectively. The process of differential equation solution can be seen as one of structure node (detail) exposition. In the numerical analysis field, the node number a large-sized element contains could be adjusted respectively in various manners by different methods. Those approaches can be classified into two categories，one is the discretized model method, featured with split node shape functions, such as the traditional FEM, the multigrid FEM, the adaptive refinement FEM, etc., another is the integrated model approach, characterized by full node shape functions, such as the wavelet FEM (WFEM), the traditional meshfree method (MFM), the traditional natural element method (NEM), isogeometric analysis method (IGAM) and the proposed multiresolution element method (MEM) etc. FEM applies the scheme of meshing and re-meshing, which is mainly relied on the empiricism, to adjust the element node number in a rough way, thus performing an irrational MRA; WFEM adopts the technique of cubic B-spline function tensor product to form the full node shape functions that are complicated to be numerically integrated and to be utilized to treat boundary conditions. MFM and NEM employ the strategy of prior artificial-selected element node layout which is also largely dependent on the empiricism. IGAM has some pitfalls like WFEM. In a word, all those above or other methods are short of the parameter-resolution level (RL) with a clear mathematical sense that can be easily used to fully alter total element node number and locate element node because they do not have a simple, clear and solid mathematical basis. However, MEM has such a simple, clear and rigorous mathematical basis that brings about the parameter RL to freely adjust total node number and locate nodes within the element. Hence, it can be said that WFEM, MFM, NEM, IGAM etc are the intermediate products in the transition of the traditional FEM from the monoresolution (discretized model) to the multiresolution (integrated model) and MEM consolidates all these irrational MRA approaches.

## 7. Conclusion and Prospective

A new multiresolution triangular plate-bending element method that has both high power of resolution and strong flexibility of analysis clarity is introduced into the field of numerical analysis. The method possesses such prominent features as follows:

1. An innovative split-full node notion is presented and a novel technique is proposed to construct a simple and clear basic full node shape function for a triangular plate-bending element, which unveils the secrets behind assembling artificially of node-related items in global matrix formation by the conventional FEM.
2. A mathematical basis for the MRA framework, that is the displacement subspace sequence, is constituted out of the scaled and shifted version of the basic full node shape function, which brings about the rational MRA concept together with the RL that defines a node uniform distribution pattern.
3. The traditional 3-node triangular plate-bending element method is a monoresolution one and also a special case of the proposed. An element of superior analysis clarity surely contains more nodes when compared with that of the inferior.

4. The RL adjusting for the multiresolution triangular plate-bending element model is laid on the rigorous mathematical basis while the meshing or remeshing for the monoresolution is based on the empiricism. Thus, the proposed element method is an rational MRA approach and can consolidate all corresponding irrational MRA approaches. As a result, selection of an appropriate RL is a time-saving task while selection of an appropriate meshing scheme time-consuming. The accuracy of a plate structure analysis is replaced by the clarity, the irrational MRA by the rational, the mesh by the RL that is the discretized model by the integrated.
5. The structural analysis clarity is actually determined by the RL, not by the mesh.
6. With advent of the new finite element method [17,18,19], the rational MRA will find a wide application in numerical solution of engineering problems in a real sense.

The upcoming work will be focused on the treatment of interface between multiresolution elements of different RL. The interface may be extended to the bridging domain in which a transitional element could be used just as PS images of different RL. The transitional element could also be constructed by the technique of scaling and shifting of the basic full node shape function to virtual or real nodes.

## References

[1] Zienkiewicz,O.C. and Taylor,R.L. The Finite Element Method. seventh ed., Butterworth-Heihemann, London. (2013)

[2] Akin, J.E. Finite Element Analysis with Error Estimators–An introduction to the FEM and adaptive error analysis for engineering student. Elseiver,Butterworth-Heihemann, Amsterdam (2005)

[3] Bartels,S**.**, Carstensen**,**C**.** A convergent adaptive finite element method for an optimal design. Numerische Mathematik, 108(3): 359-385 (2008).

[4] Rivara,M.C. Local modification of meshes for adaptive and/or multigrid finite-element methods. Journal of Computational and Applied Mathematics, 36(1): 79-89 (1991)

[5] Shaidueov,V.V. Multigrid Methods for Finite Elements. Kluwer, Amsterdam(1995)

[6] Xiang, J.W., Chen, X.F., He, Y.M and He, Z.J. The construction of plane elastomechanics and Mindlin plate elements of B-spline wavelet on the interval. Finite Elements in Analysis and Design, **42**: 1269-1280 (2006)

[7] He, Z.J., Chen, X.F., and Li, B. Theory and engineering application of wavelet finite element method, Science Press, Beijing (2006)

[8] Yu, Y., Lin, Q. Y. and Yang, C. A 3D shell-like approach using element-free Galerkin method for analysis of thin and thick plate structures. Acta Mechanica Sinica, **29**: 85-98 (2013)

[9] Liu, H.S. and Fu, M.W. Adaptive reproducing kernel particle method using gradient indicator for elasto-plastic deformation. Engineering Analysis with Boundary Elements, **37**: 280–292 (2013)

[10] Sukumar, N., Moran, B. and Belytschko, T. The natural elements method in solid mechanics. International Journal of Numerical Methods in Engineering.**43**: 839-887 (1998)

[11] Sukumar, N., Moran, B. and Semenov, A.Y. Natural neighbor Galerkin methods, International Journal of Numerical Methods in Engineering. **50**: 1-27(2001)

[12] H. Lian, P. Kerfriden, S.P.A. Bordas , Shape optimization directly from CAD: An isogeometric boundary element approach using T-splines, Computer Methods in Applied Mechanics and. Engineering. **317**: 1–41(2017).

[13] T.J.R.HughesJ.A.CottrellY.Bazilevs**,** Isogeometric analysis: CAD, finite elements, NURBS, exact geometry and mesh refinement. Computer Methods in Applied Mechanics and. Engineering. **194(**39-41): 4135–4195(2005) .

[14] Long, Y.Q., Cen, S. and Long, Z.F. Advanced Finite Element Method in Structural Engineering. Springer-Verlag GmbH: Berlin, Heidelberg; Tsinghua University Press: Beijing, 2009

[15] Bergan P G, Felippa C A. A triangular membrane element with rotational degrees of freedom. Computer Methods in Applied Mechanics and Engineering,,**50**:25-69,(1985)

[16] Wang H.X., Cen S. Locking-free Thick Plate Bending Triangle Element. Journal of Nantong University, 2(2):5-11 (2003) (in Chinese)

[17] Xia, Y. M., Liu, Y. X., Chen, S.L.and Tan, G. A rectangular shell element formulation with a new multi- resolution analysis. Acta Mechanica Solida Sinca. **27(6)**:612-625(2014)

[18] Xia, Y. M. A multiresolution finite element method based on a new quadrilateral plate element. Journal of Coupled Systems and Multiscale Dynamics. **2(2)**:52-61 (2014)

[19] Xia, Y. M., Chen, S.L. A hexahedron element formulation with a new multi-resolution analysis. Science China (Physics, Astronomy and Mechanics), **58(1)**: 014601-10 (2015)